\documentclass[english]{article}
\usepackage[T1]{fontenc}
\usepackage[latin9]{inputenc}
\usepackage{geometry}
\usepackage{mathrsfs}
\usepackage{amsmath}
\usepackage{amssymb}
\usepackage{stackrel}
\usepackage{setspace}
\PassOptionsToPackage{normalem}{ulem}
\usepackage{ulem}
\makeatletter
\@ifundefined{date}{}{\date{}}
\makeatother

\usepackage{babel}
\begin{document}

\title{Furstenberg's Times 2, Times 3 Conjecture (a Short Survey)}

\author{Matan Tal\\
The Hebrew University of Jerusalem}
\maketitle
\begin{abstract}
The following is a concise exposition of the conjecture and three
of its proofs for the case of positive entropy by D. Rudolph \cite{key-20}
, B. Host \cite{key-15} and W. Parry \cite{key-19}. A simpler theorem
of R. Lyons \cite{key-17} - preceding them - is also presented and
proved. This is a survey, no new results are introduced.
\end{abstract}

\section{Introduction}

Throughout this exposition, $\mathbb{T}$ will denote for us $\mathbb{R}/\mathbb{Z}$.
Let $p,q>1$ be multiplicatively independent integers, i.e. they are
not both powers of the same integer. H. Furstenberg proved in \cite{key-11}
(1967) the following theorem:\\

\textbf{Theorem 1.1:} If $F\subseteq\mathbb{T}$ is infinite, closed
and invariant under multiplication by both $p$ and $q$ then $F=\mathbb{T}$.\\

An equivalent formulation of the theorem is that the orbit of every
irrational point under the semi-group $\left\langle p,q\right\rangle $
is dense. Furstenberg also raised the following conjecture which is
the measurable analog of Theorem 1.1:\\

\textbf{Conjecture 1.2 (Furstenberg's Times 2, Times 3 Conjecture):}
If $\mu\in\text{\ensuremath{M\left(\mathbb{T}\right)}}$ (the probability
measures on $\left(\mathbb{T},\mathcal{B}_{\mathbb{T}}\right)$) is
an atomless invariant measure under multiplication by both $p$ and
$q$ then $\mu=\lambda$ (Lebesgue measure).\\

(Atomic ergodic measures which are likewise invariant of course do
exist - those are the uniform measures on the orbits of the action
of the multiplicative semi-group $\left\langle p,q\right\rangle $
on any rational number with a denominator prime to both $p$ and $q$.
An empirical study concerning such measures can be found in \cite{key-7}.)\\

By deploying ergodic decomposition one sees that requiring $\mu$
to be also ergodic with respect to the action of the multiplicative
semi-group $\left\langle p,q\right\rangle $ results in an equivalent
conjecture. It is tempting to think the conjecture immediately implies
Theorem 1.1, by taking an ergodic invariant measure supported on the
orbit closure of an irrational point. However, In order to prove that
there exists such a measure that is not supported on a finite set
of rationals, it seems one must argue in a similar fashion as the
proof of Theorem 1.1 itself.\\

Contemporary knowledge about the conjecture's validness is more or
less summarized by Rudolph's Theorem \cite{key-20}.\\

\textbf{Theorem 1.3 (Rudolph's Theorem):} Let $p,q>1$ be relatively
prime integers. If $\mu\in\text{\ensuremath{M\left(\mathbb{T}\right)}}$
is an invariant measure under multiplication by both $p$ and $q$,
ergodic under the action of the multiplicative semi-group $\left\langle p,q\right\rangle $,
and, moreover, there exists $r\in\left\langle p,q\right\rangle $
with $h_{\mu}\left(r\right)>0$ then $\mu=\lambda$.\\

We shall see (Proposition 3.1) that there exists such an $r$ if and
only if multiplication by every element of $\left\langle p,q\right\rangle $
has positive entropy. Therefore one could have equally required $h_{\mu}\left(p\right)>0$
instead.\\

Actually, A. Johnson, Rudolph's student, improved Rudolph's proof
to apply for multiplicatively independent integers $p,q$ \cite{key-16}
(this full version of the theorem is called Rudolph-Johnson Theorem),
but since Lyon's and Host's ideas presented here apply only to relatively
prime integers $p,q$ we stated the theorem as we did for our convenience.
(Notice that by passing to the invertible extension, Rudolph's theorem
for the relatively prime integers $p,q>1$ implies its validity also
for $p^{m_{1}}q^{n_{1}},p^{m_{2}}q^{n_{2}}$ where $m_{1},m_{2},n_{1},n_{2}\in\mathbb{N}$
and $m_{1}n_{2}-n_{1}m_{2}\neq0$.)\\

Not much is known about the case of zero entropy - in \cite{key-1},
the existence of such an invariant measure without atoms is shown
to be equivalent to the existence of two partitions that satisfy some
condition that is formulated using only the Lebesgue measure of the
circle. Additional known facts about the conjecture (and Theorem 1.1)
that are not included in this exposition include its generalizations
to higher dimensional tori and other connected compact Abelian groups
\cite{key-4-1,key-5,key-8}, some quantitative estimates \cite{key-6},
treatments of related conjectures and other proofs of Rudolph's Theorem
that exist such as that of Feldman \cite{key-4-3} (a weaker result
than Rudolph-Johnson Theorem) and that of Hochman and Shmerkin \cite{key-13,key-14}
(a proof of Rudolph-Johnson Theorem - it uses ergodic theoretical
methods associated with fractal geometry originating from \cite{key-12}).
We shall also not elaborate on the importance of the subject, but
only mention here that apart from raising a very natural question,
the conjecture also serves as the leading toy example for the more
general phenomenon of measure rigidity of higher rank hyperbolic actions
(for which contemporary understanding is in a similar situation \cite{key-22,key-2,key-3,key-4}).\\

Apart from the application of Prop 3.1 in section 4, the following
sections are independent from one another.\\

I wish to thank my doctoral advisors Prof. H. Furstenberg and Prof.
T. Meyerovitch for having the patience to hear my presentation of
most of this survey and for some valuable remarks.

\section{Lyon's Theorem}

The first significant advancement towards a proof of the Times 2,
Times 3 conjecture was made by R. Lyons in \cite{key-17} (1988).\\

\textbf{Definiton:} A measure preserving system $\left(X,T,\mu\right)$
is $K$-mixing (or $T$-exact) if for every

$g\in L^{2}\left(X,\mu\right)$ the limit $\left\langle f\circ T^{n},g\right\rangle -\left\langle f,1\right\rangle \left\langle 1,g\right\rangle \rightarrow0$
as $n\rightarrow\infty$ exists uniformly for $\left\Vert f\right\Vert _{2}\leq1$.\\

\textbf{Theorem 2.1:} Let $p,q>1$ be relatively prime integers. If
$\delta_{0}\neq\mu\in\text{\ensuremath{M\left(\mathbb{T}\right)}}$
is an invariant measure under multiplication by both $p$ and $q$
and $K$-mixing for the times $p$ transformation then $\mu=\lambda$.\\

So let $p,q>1$ be relatively prime integers. Lyons proves two elementary
number theoretical lemmas.\\

\textbf{Lemma 2.2:} There exist $A,d,L\in\mathbb{N}$ such that $q^{A}=dp^{L}+1$
and $L\geq2$, $\gcd\left(d,p\right)=1$.\\

\textbf{Lemma 2.3:} Assuming the minimal $A$ for which $d,L$ exist
satisifying Lemma 2.2, then for all $l\geq L$ the order of $q$ modulo
$p^{l}$ is $p^{l-L}A$.\\

So $q^{x}$ attains $p^{l-L}A$ values Modulo $p^{l}$, and this implies
a conclusion.\\

\textbf{Conclusion 2.4:} $q^{x}\equiv b\,\left(\mod p^{l}\right)$
has a solution if and only if $q^{x}\equiv b\,\left(\mod p^{L}\right)$
does.\\

The proof of Theorem 2.1 proceeds as follows. From Conclusion 2.4
we know that there exists a strictly increasing sequnce $n_{j}$ for
$j\geq L$ such that $q^{n_{j}}\equiv p^{j}+1\,\left(\mod p^{2j}\right)$
(since

$q^{x}\equiv p^{j}+1\,\left(\mod p^{L}\right)$ has a solution for
every such $j$), i.e. $q^{n_{j}}=d_{j}p^{2j}+p^{j}+1$ for some integer
$d_{j}$.\\

For $0\neq m\in\mathbb{Z}$, we need to prove that $\hat{\mu}\left(m\right)=0$.
Writing\\

$\hat{\mu}\left(m\right)=\hat{\mu}\left(mq^{n_{j}}\right)=\hat{\mu}\left(m\left(d_{j}p^{2j}+p^{j}+1\right)\right)=\hat{\mu}\left(\left(md_{j}p^{j}+m\right)p^{j}+m\right)$,\\

and substituting twice two characters for the functions $f,g$ in
the definition of $K$-mixing, we obtain:\\

$\hat{\mu}\left(m\right)=\hat{\mu}\left(m\right)\lim_{j\rightarrow\infty}\hat{\mu}\left(md_{j}p^{j}+m\right)=\hat{\mu}\left(m\right)^{2}\lim_{j\rightarrow\infty}\hat{\mu}\left(md_{j}\right)$.\\

Assuming to the contrary that $\hat{\mu}\left(m\right)\neq0$, the
conclusion $\left|\hat{\mu}\left(m\right)\right|=1$ is forced upon
us. This means $\int_{\mathbb{T}}e^{2m\pi xi}\,d\mu\left(x\right)$
has no cancellations whatsoever, i.e. $\mu$ is supported on finitely
many points. This is a contradiction since then our system is not
even mixing. 

\section{Rudolph's Theorem}

Motivated by the implicit role positive entropy plays in Lyons' result
\footnote{The $K$-mixing hypothesis in Theorem 2.1, implies triviality of the
Pinsker $\sigma$-algebra of $\sigma_{p}$, the multiplication by
$p$ transformation. Indeed, If a set $C$ belongs to the Pinsker
$\sigma$-algebra of $\sigma_{p}$, then - modulo $\mu$ - it is a
tail event of the partition $\left\{ \left[\frac{j}{p},\frac{j+1}{p}\right)\right\} _{j=0}^{p-1}$.
This means that for every $n>0$ there exists a Borel measurable set
$C_{n}$ such that $C=\sigma_{p}^{-n}\left(C_{n}\right)$. Taking
$f_{n}=1_{C_{n}}$ and $g=1_{C}$ and using $K$-mixing we deduce
$\mu\left(C\right)$ equals either $0$ or $1$.}, D. Rudolph published his theorem (Theorem 1.3) two years later in
1990 \cite{key-20}. In this section we present his proof (with very
slight modifications). \\

We begin with a fundamental observation by Rudolph.\\

\textbf{Proposition 3.1:} Let $a,b\in\mathbb{N}$ and $\sigma_{a}$,
$\sigma_{b}$ the maps of multiplication of $\mathbb{T}$ by each
of them respectively. If $\mu\in M\left(\mathbb{T}\right)$ is an
invariant measure under $\sigma_{a}$ and $\sigma_{b}$, then $\log b\cdot h_{\mu}\left(\sigma_{a}\right)=\log a\cdot h_{\mu}\left(\sigma_{b}\right)$.\\

\textbf{Proof:} Assume $a,b>1$ for otherwise there is nothing to
prove. Denote by $\xi$ the partition $\left\{ \left[\frac{j}{ab},\frac{j+1}{ab}\right)\right\} _{j=0}^{ab-1}$.
The lengths of the intervals composing the partition $\vee_{i=0}^{l-1}\sigma_{a}^{-i}\left(\xi\right)$
is $\frac{1}{a^{l}b}$, and of $\vee_{i=0}^{m-1}\sigma_{b}^{-i}\left(\xi\right)$
is $\frac{1}{ab^{m}}$. Choosing $m\left(l\right)=\left\lfloor \log_{b}a^{l-1}\right\rfloor =\left\lfloor \frac{\left(l-1\right)\log a}{\log b}\right\rfloor $
the length of the intervals composing $\vee_{i=0}^{m\left(l\right)-1}\sigma_{b}^{-i}\left(\xi\right)$
is between $\frac{1}{a^{l}b}$ and $\frac{1}{a^{l}}$, and hence none
of which can intersect more than $b+1$ of the intervals composing
$\vee_{i=0}^{l-1}\sigma_{a}^{-i}\left(\xi\right)$. Thus\\

$H_{\mu}\left(\vee_{i=0}^{l-1}\sigma_{a}^{-i}\left(\xi\right)\right)\leq H_{\mu}\left(\left(\vee_{i=0}^{l-1}\sigma_{a}^{-i}\left(\xi\right)\right)\vee\left(\vee_{i=0}^{m\left(l\right)-1}\sigma_{b}^{-i}\left(\xi\right)\right)\right)$\\
$=H_{\mu}\left(\vee_{i=0}^{m\left(l\right)-1}\sigma_{b}^{-i}\left(\xi\right)\right)+H_{\mu}\left(\vee_{i=0}^{l-1}\sigma_{a}^{-i}\left(\xi\right)\,|\,\vee_{i=0}^{m\left(l\right)-1}\sigma_{b}^{-i}\left(\xi\right)\right)\leq H_{\mu}\left(\vee_{i=0}^{m\left(l\right)-1}\sigma_{b}^{-i}\left(\xi\right)\right)+\log\left(b+1\right)$.\\

By Dividing both sides by $l$ we get

$\frac{1}{l}H_{\mu}\left(\vee_{i=0}^{l-1}\sigma_{a}^{-i}\left(\xi\right)\right)\leq\frac{m\left(l\right)}{l}\frac{1}{m\left(l\right)}H_{\mu}\left(\vee_{i=0}^{m\left(l\right)-1}\sigma_{b}^{-i}\left(\xi\right)\right)+\frac{1}{l}\log\left(b+1\right)$,
and by letting $l\rightarrow\infty$ we conclude $\log b\cdot h_{\mu}\left(\sigma_{a}\right)\leq\log a\cdot h_{\mu}\left(\sigma_{b}\right)$.
The reverse inequality is obtained similarly. $\blacksquare$\\

An immediate conclusion is that if $\mu\in M\left(\mathbb{T}\right)$
is invariant under a semi-group of natural numbers different than
$1$, than the entropy is either positive for all its elements or
zero for all its elements. \\

Now let $p,q,\mu$ be as in Rudolph's Theorem (Theorem 1.3), $T_{0}:\mathbb{T}\rightarrow\mathbb{T}$
the multiplication by $p$ and $S_{0}:\mathbb{T}\rightarrow\mathbb{T}$
by $q$. \\

\subsection{Preparations}

We construct a subshift of finite type (SFT) $X$ for $\mathbb{Z}^{2}$
with the alphabet $\Lambda=\left\{ 0,1,2,\dots,pq-1\right\} $. Denote
by $T$ the action of $\left(1,0\right)$ and by $S$ that of $\left(0,1\right)$.
The rule of the subshift allows $i\in\Lambda$ to horizontally precede
$j\in\Lambda$ if $\left(\frac{i}{pq},\frac{i+1}{pq}\right)\cap T_{0}^{-1}\left(\left(\frac{j}{pq},\frac{j+1}{pq}\right)\right)\neq\emptyset$
(The pre-image on the right is composed of $p$ intervals of length
$\frac{1}{p^{2}q}$), and likewise $i\in\Lambda$ to vertically precede
$j\in\Lambda$ if $\left(\frac{i}{pq},\frac{i+1}{pq}\right)\cap S_{0}^{-1}\left(\left(\frac{j}{pq},\frac{j+1}{pq}\right)\right)\neq\emptyset$.
The purpose of this subsection is to clear the ground and understand
properties of $X$ and its exact relation to our original system on
the circle. The strategy of the proof is introduced only in the next
subsection.\\

Associating to this the horizontal adjacency matrix $M_{T}$ and the
vertical one $M_{S}$ we deduce that $M_{T}M_{S}=M_{S}M_{T}$ is just
the $pq\times pq$ matrix with all entries equal to $1$ (this is
because the diagonal of the subshift is the full one-dimensional shift
$\Lambda$). Notice then while we have complete freedom to choose
the values of a positively infinite diagonal ray (of $\left(1,1\right)$
differences), it then determines all values of its corresponding quadrant.\\

Define $\varphi:X\rightarrow\mathbb{T}$ by $\varphi\left(x\right)=\stackrel[n=0]{\infty}{\sum}\frac{x\left(n,n\right)}{\left(pq\right)^{n+1}}$.
It is equivariant and, on the $\left(\mathbb{N}\cup\left\{ 0\right\} \right)^{2}$-restricted
version of $X$ (taking only the part of $X$ which is in the non-negative
quadrant) $\varphi$ is almost one-to-one. The bad set (which has
two-point fibers) is $V=\left\{ \frac{t}{\left(pq\right)^{n}}\,:\,0\leq t,n\in\mathbb{Z}\right\} $.
In particular, this implies that, in $X\setminus\varphi^{-1}\left(V\right)$,
the values on any positively infinite horizontal ray (of $\left(1,0\right)$
differences) determine all values above it, and the values on any
positively infinite vertical ray (of $\left(0,1\right)$ differences)
determine all values to its right.\\

Assuming $\mu$ is not supported on $V$ (that is $\mu\neq\delta_{0}$)
we can lift it up to a unique $\left\langle T,S\right\rangle $\textendash ergodic
measure on $X\setminus\varphi^{-1}\left(V\right)$ which we will denote
by $\tilde{\mu}$. Denote by $\alpha$ the partition of $X$ by the
value at the coordinate $\left(0,0\right)$.\\

\textbf{Lemma 3.2: }$h_{\tilde{\mu}}\left(T\,|\,\mathcal{A}\right)=h_{\tilde{\mu}}\left(T,\alpha\,|\,\mathcal{A}\right),\,h_{\tilde{\mu}}\left(S\,|\,\mathcal{A}\right)=h_{\tilde{\mu}}\left(S,\alpha\,|\,\mathcal{A}\right)$
where $\mathcal{A}$ is any sub-$\sigma$-algebra invariant under
both $T$ and $S$ (i.e. $S^{-1}\mathcal{A}=T^{-1}\mathcal{A}=\mathcal{A}$).
\\

\textbf{Proof:} Just notice that $S^{m}\stackrel[i=-\infty]{\infty}{\vee}T^{-i}\left(\alpha\right)\,\underset{m\rightarrow\infty}{\nearrow}$
the $\sigma$-algebra of $X$, so by Kolmogorov-Sinai Theorem $h_{\tilde{\mu}}\left(T\,|\,\mathcal{A}\right)=\lim_{m\rightarrow\infty}h_{\tilde{\mu}}\left(T,S^{m}\left(\alpha\right)\,|\,\mathcal{A}\right)=h_{\tilde{\mu}}\left(T,\alpha\,|\,\mathcal{A}\right)$.
The proof that $h_{\tilde{\mu}}\left(S\,|\,\mathcal{A}\right)=h_{\tilde{\mu}}\left(S,\alpha\,|\,\mathcal{A}\right)$
is similar. $\blacksquare$\\

The fact that $\gcd\left(p,q\right)=1$ is taken advantage of only
through the following observation and in the proof of Lemma 3.6.\\

\textbf{Lemma 3.3:} For any values $x\left(n,k\right)$ on all $n\geq1,\,k\geq0$
together with $x\left(0,m_{0}\right)$ for some $m_{0}\geq0$ all
satisfying the SFT rule with their neighbors among those coordintates,
there exists $x$s in $X$ with those values and they are all equal
also on $\left(0,k\right)$ for all $k\geq0$. Similarly, for any
values $x\left(n,k\right)$ on all $n\geq0,\,k\geq1$ together with
$x\left(m_{0},0\right)$ for some $m_{0}\geq0$ all satisfying the
SFT rule with their neighbors among those coordintates, there exists
$x$s in $X$ with those values and they are all equal also on $\left(n,0\right)$
for all $n\geq0$.\\

\textbf{Proof:} We only prove the first statement (the proof of the
second one is similar). $x\left(0,m_{0}\right)$ determines $x\left(0,m_{0}+1\right)$
since $M_{T}M_{S}$ has all its entries equal to $1$. It then suffices
to prove that $x\left(0,0\right)$ is determined when $m_{0}=1$.\\

$T_{0}^{-1}\left(\varphi\left(T\left(x\right)\right)\right)=\left\{ \frac{\varphi\left(T\left(x\right)\right)+i}{p}\,:\,0\leq i<p\right\} $
and thus $S_{0}\left(T_{0}^{-1}\left(\varphi\left(T\left(x\right)\right)\right)\right)=\left\{ \frac{\varphi\left(T\left(x\right)\right)q+iq}{p}\,:\,0\leq i<p\right\} $.
Since $q$ is invertible in the ring $\mathbb{Z}/p\mathbb{Z}$ there
exists a one-to-one correspondance between allowed entries for $x\left(0,0\right)$
and allowed entries for $x\left(0,1\right)$. Thus the choice of the
former determines the latter. $\blacksquare$\\

\textbf{Lemma 3.4:} For $\tilde{\mu}$-a.e. $x$ the events $\left\{ x'\in X\,:\,x'\left(-1,0\right)=k_{1}\right\} $
and $\left\{ x'\in X\,:\,x'\left(0,-1\right)=k_{2}\right\} $ for
$0\leq k_{1},k_{2}<pq$ are independent when conditioned on the first
quadrant 

(i.e. the quadrant $\stackrel[i=0]{\infty}{\vee}\stackrel[j=0]{\infty}{\vee}T^{-i}S^{-j}\left(\alpha\right)=\stackrel[i=0]{\infty}{\vee}T^{-i}S^{-i}\left(\alpha\right)$).\\

\textbf{Proof:} It suffices to prove that

$H_{\tilde{\mu}}\left(T^{-1}\left(\alpha\right)\vee S^{-1}\left(\alpha\right)\,|\,\stackrel[i=1]{\infty}{\vee}T^{-i}S^{-i}\left(\alpha\right)\right)=H_{\tilde{\mu}}\left(T^{-1}\left(\alpha\right)\,|\,\stackrel[i=1]{\infty}{\vee}T^{-i}S^{-i}\left(\alpha\right)\right)+H_{\tilde{\mu}}\left(S^{-1}\left(\alpha\right)\,|\,\stackrel[i=1]{\infty}{\vee}T^{-i}S^{-i}\left(\alpha\right)\right)$.
Since a positive diagonal ray determines a quadrant and also by applying
Lemma 3.3:

$\alpha\vee\left(\stackrel[i=1]{\infty}{\vee}T^{-i}S^{-i}\left(\alpha\right)\right)=T^{-1}\left(\alpha\right)\vee S^{-1}\left(\alpha\right)\vee\left(\stackrel[i=1]{\infty}{\vee}T^{-i}S^{-i}\left(\alpha\right)\right)$.
Hence 

$H_{\tilde{\mu}}\left(\alpha\,|\,\stackrel[i=1]{\infty}{\vee}T^{-i}S^{-i}\left(\alpha\right)\right)=H_{\tilde{\mu}}\left(T^{-1}\alpha\vee S^{-1}\alpha\,|\,\stackrel[i=1]{\infty}{\vee}T^{-i}S^{-i}\left(\alpha\right)\right)$

$\leq H_{\tilde{\mu}}\left(T^{-1}\left(\alpha\right)\,|\,\stackrel[i=1]{\infty}{\vee}T^{-i}S^{-i}\left(\alpha\right)\right)+H_{\tilde{\mu}}\left(S^{-1}\left(\alpha\right)\,|\,\stackrel[i=1]{\infty}{\vee}T^{-i}S^{-i}\left(\alpha\right)\right)=h_{\tilde{\mu}}\left(T\right)+h_{\tilde{\mu}}\left(S\right)$
(the last equality is by Lemma 3.2). So we want to show the inequality
is in fact an equality, namely that $H_{\tilde{\mu}}\left(\alpha\,|\,\stackrel[i=1]{\infty}{\vee}T^{-i}S^{-i}\left(\alpha\right)\right)=h_{\tilde{\mu}}\left(T\right)+h_{\tilde{\mu}}\left(S\right)$.
This is an immediate consequence of Prop. 3.1 (because $H_{\tilde{\mu}}\left(\alpha\,|\,\stackrel[i=1]{\infty}{\vee}T^{-i}S^{-i}\left(\alpha\right)\right)=h_{\tilde{\mu}}\left(TS\right)$),
but let us also present a computation that does not make use of it.\\

$H_{\tilde{\mu}}\left(\alpha\,|\,\stackrel[i=1]{\infty}{\vee}T^{-i}S^{-i}\left(\alpha\right)\right)=H_{\tilde{\mu}}\left(\alpha\vee\left(\stackrel[i=1]{\infty}{\vee}T^{-i}\left(\alpha\right)\right)\,|\,\stackrel[i=1]{\infty}{\vee}T^{-i}S^{-i}\left(\alpha\right)\right)$
since the diagonal ray 

$\stackrel[i=0]{\infty}{\vee}T^{-i}S^{-i}\left(\alpha\right)=\alpha\vee\left(\stackrel[i=1]{\infty}{\vee}T^{-i}S^{-i}\left(\alpha\right)\right)$
determines the corresponding quadrant which includes $\stackrel[i=1]{\infty}{\vee}T^{-i}\alpha$.
But

$H_{\tilde{\mu}}\left(\alpha\vee\left(\stackrel[i=1]{\infty}{\vee}T^{-i}\left(\alpha\right)\right)\,|\,\stackrel[i=1]{\infty}{\vee}T^{-i}S^{-i}\left(\alpha\right)\right)=H_{\tilde{\mu}}\left(\alpha\,|\,\stackrel[i=1]{\infty}{\vee}\stackrel[j=0]{\infty}{\vee}T^{-i}S^{-j}\left(\alpha\right)\right)+H_{\tilde{\mu}}\left(T^{-1}\left(\alpha\right)\,|\,\stackrel[i=1]{\infty}{\vee}T^{-i}S^{-i}\left(\alpha\right)\right)$
(because $\left(T^{-1}\alpha\right)\vee\left(\stackrel[i=1]{\infty}{\vee}T^{-i}S^{-i}\left(\alpha\right)\right)=\stackrel[i=1]{\infty}{\vee}\stackrel[j=0]{\infty}{\vee}T^{-i}S^{-j}\left(\alpha\right)$)
and
\[
H_{\tilde{\mu}}\left(\alpha\,|\,\stackrel[i=1]{\infty}{\vee}\stackrel[j=0]{\infty}{\vee}T^{-i}S^{-j}\left(\alpha\right)\right)+H_{\tilde{\mu}}\left(T^{-1}\left(\alpha\right)\,|\,\stackrel[i=1]{\infty}{\vee}T^{-i}S^{-i}\left(\alpha\right)\right)
\]

\[
=H_{\tilde{\mu}}\left(\alpha\,|\,\stackrel[i=1]{\infty}{\vee}T^{-i}\left(\alpha\right)\right)+H_{\tilde{\mu}}\left(\alpha\,|\,\stackrel[i=1]{\infty}{\vee}S^{-i}\left(\alpha\right)\right)=h_{\tilde{\mu}}\left(T\right)+h_{\tilde{\mu}}\left(S\right).\,\blacksquare
\]
\\

\textbf{Lemma 3.5:} If $\mathcal{A}$ is a sub-$\sigma$-algebra invariant
under both $T$ and $S$ (i.e. $T^{-1}\mathcal{A}=S^{-1}\mathcal{A}=\mathcal{A}$)
then 

$h_{\tilde{\mu}}\left(T;\mathcal{A}\right)=\frac{\log p}{\log q}h_{\tilde{\mu}}\left(S;\mathscr{\mathcal{A}}\right)$
- where $h_{\tilde{\mu}}\left(S;\mathscr{\mathcal{A}}\right)$ denotes
the entropy of the action of $S$ on $\left(X,\mathscr{\mathcal{A}},\tilde{\mu}\right)$.\\

\textbf{Proof:} By lifting Prop. 3.1 from the circle to $X$ (through
$\varphi$) one concludes that $h_{\tilde{\mu}}\left(T\right)=\frac{\log p}{\log q}h_{\tilde{\mu}}\left(S\right)$,
and by a similar proof to Prop. 3.1 that $h_{\tilde{\mu}}\left(T\,|\,\mathcal{A}\right)=\frac{\log p}{\log q}h_{\tilde{\mu}}\left(S\,|\,\mathscr{\mathcal{A}}\right)$.
The Abramov-Rokhlin formula yields the desired result. $\blacksquare$\\

\subsection{The Heart of the Proof}

Supposing $\mu$ is not Lebesgue measure $\lambda$, we want to show
that $h_{\tilde{\mu}}\left(T\right)=0$. The strategy of the proof
is based on the insight that it is enough to find a sub-$\sigma$-algebra
$\mathcal{A}$ invariant under both $T$ and $S$ (i.e. $T^{-1}\mathcal{A}=S^{-1}\mathcal{A}=\mathcal{A}$)
such that (i) $h_{\tilde{\mu}}\left(S;\mathscr{\mathcal{A}}\right)=0$
and (ii) $\alpha\underset{\mod\tilde{\mu}}{\subseteq}\mathcal{A}\vee\left(\stackrel[i=1]{\infty}{\vee}T^{-i}\left(\alpha\right)\right)$.
This is because then $h_{\tilde{\mu}}\left(T\right)=h_{\tilde{\mu}}\left(T\,|\,\mathcal{A}\right)+h_{\tilde{\mu}}\left(T;\mathcal{A}\right)=h_{\tilde{\mu}}\left(T;\mathcal{A}\right)$
(the first equality is by the Rokhlin-Abramov formula, and the second
by Lemma 3.2 and (ii)). But $h_{\tilde{\mu}}\left(T;\mathcal{A}\right)=\frac{\log p}{\log q}h_{\tilde{\mu}}\left(S;\mathscr{\mathcal{A}}\right)=0$
by lemma 3.5. The rest of this subsection is dedicated to finding
such an $\mathcal{A}$ and proving that it satisfies those two desired
properties.\\

For this purpose, given any fixed $n\in\mathbb{N}$ we define a function
$\nu_{n}\left(\cdot\right):X\rightarrow M\left(\mathbb{T}\right)$
for $\tilde{\mu}$-a.e. $x$ by the formula $\nu_{n}\left(x\right)=\stackrel[t=0]{p^{n}-1}{\sum}a_{n}^{\left(t\right)}\left(x\right)\cdot\delta_{\frac{t}{p^{n}}}$
where

$a_{n}^{\left(t\right)}\left(x\right)=\tilde{\mu}_{x}^{\stackrel[i=0]{\infty}{\vee}T^{-i}\alpha}\left(\left\{ x'\in X\,:\,x'(-j,0)=T^{n}\circ\varphi^{-1}\left(\varphi\left(T^{-n}\left(x\right)\right)+\frac{t}{p^{n}}\right)\left(-j,0\right)\,for\,all\,1\leq j\leq n\right\} \right)$
($\tilde{\mu}_{x}^{\stackrel[i=0]{\infty}{\vee}T^{-i}\left(\alpha\right)}$
is the conditional measure at $x$ with respect to $\stackrel[i=0]{\infty}{\vee}T^{-i}\left(\alpha\right)$).
$\nu_{n}$ is measurable according to the weak-{*} topology on its
range (equivalently, all $a_{n}^{\left(t\right)}\left(x\right)$ are
real measurable functions). \\

\textbf{Lemma 3.6:} $a_{n}^{\left(qt\,\mod p^{n}\right)}\left(S\left(x\right)\right)=a_{n}^{\left(t\right)}\left(x\right)$
for $\tilde{\mu}$-a.e. $x$.\\

\textbf{Proof:} We will explain only the case of $n=1$ as the general
case follows by a similar inductive argument \footnote{Using also the trivial fact that if events $A_{1},A_{2},A_{3},B_{1},B_{2},B_{3}$
in a probability space satisfy

$\mathbb{P}\left(A_{2}\,|\,A_{3}\right)=\mathbb{P}\left(B_{2}\,|\,B_{3}\right)$
and $\mathbb{P}\left(A_{1}\,|\,A_{2}\cap A_{3}\right)=\mathbb{P}\left(B_{1}\,|\,B_{2}\cap B_{3}\right)$
then $\mathbb{P}\left(A_{1}\cap A_{2}\,|\,A_{3}\right)=\mathbb{P}\left(B_{1}\cap B_{2}\,|\,B_{3}\right)$.}. So why does $a_{1}^{\left(qt\,\mod p\right)}\left(S\left(x\right)\right)=a_{1}^{\left(t\right)}\left(x\right)$?
\\

$a_{1}^{\left(qt\,\mod p\right)}\left(S\left(x\right)\right)=\tilde{\mu}_{S\left(x\right)}^{\stackrel[i=0]{\infty}{\vee}T^{-i}\left(\alpha\right)}\left(\left\{ x'\in X\,:\,x'(-1,0)=T\circ\varphi^{-1}\left(\varphi\left(T^{-1}S\left(x\right)\right)+\frac{qt}{p}\right)\left(-1,0\right)\right\} \right)$\\

$=\tilde{\mu}_{x}^{\stackrel[i=0]{\infty}{\vee}S^{-1}T^{-i}\left(\alpha\right)}\left(\left\{ x'\in X\,:\,x'(-1,0)=T\circ\varphi^{-1}\left(\varphi\left(T^{-1}\left(x\right)\right)+\frac{t}{p}\right)\left(-1,0\right)\right\} \right)$\\

$=\tilde{\mu}_{x}^{\alpha\vee\left(\stackrel[i=0]{\infty}{\vee}S^{-1}T^{-i}\left(\alpha\right)\right)}\left(\left\{ x'\in X\,:\,x'(-1,0)=T\circ\varphi^{-1}\left(\varphi\left(T^{-1}\left(x\right)\right)+\frac{t}{p}\right)\left(-1,0\right)\right\} \right)$\\

where in the last equality we use Lemma 3.4. Since $\alpha\vee\left(\stackrel[i=0]{\infty}{\vee}S^{-1}T^{-i}\left(\alpha\right)\right)=\stackrel[i=0]{\infty}{\vee}T^{-i}\left(\alpha\right)$,
this is indeed $a_{1}^{\left(t\right)}\left(x\right)$. $\blacksquare$\\

We are now ready to define our invariant $\sigma$-algebra $\mathcal{A}$
as the minimal $T$-invariant (i.e.

$T^{-1}\mathcal{A}=\mathcal{A}$) sub-$\sigma$-algebra on $X$ for
which all the functions $\nu_{n}$ are measurable. By Lemma 3.6 it
is also $S$ invariant (i.e. $S^{-1}\mathcal{A}=\mathcal{A}$). Notice
that $\mathcal{A}_{k}\nearrow\mathcal{A}$ where - for every $k\geq0$
- $\mathcal{A}_{k}$ is the minimal sub-$\sigma$-algebra for which
the function $\nu_{2k+1}\circ T^{k}$ is measurable.\\

\textbf{Lemma 3.7:} There is $j_{k}$ such that $S^{j_{k}}\left(A\right)=A$
for all $A\in\mathcal{A}_{k}$.\\

\textbf{Proof:} Take $j_{k}$ to be the order of $q$ in $\left(\mathbb{Z}/p^{2k+1}\mathbb{Z}\right)^{\times}$.
By Lemma 3.6,

$a_{2k+1}^{\left(q^{j_{k}}t\,\mod p^{2k+1}\right)}\left(S^{j_{k}}T^{k}\left(x\right)\right)=a_{2k+1}^{\left(t\right)}\left(T^{k}\left(x\right)\right)$,
and so $a_{2k+1}^{\left(t\,\mod p^{2k+1}\right)}\left(S^{j_{k}}T^{k}\left(x\right)\right)=a_{2k+1}^{\left(t\right)}\left(T^{k}\left(x\right)\right)$
which means the function $\nu_{2k+1}\circ T^{k}$ stays invariant
under the action of $S^{j_{k}}$. \textbf{$\blacksquare$}\\

We can now conclude desired property (i) of $\mathscr{\mathcal{A}}$\textbf{.}\\

\textbf{Corollary 3.8:} Each ergodic component with respect to the
action of $S$ on the factor corrresponding to $\mathcal{A}$ has
a rational pure point spectrum. In particular,\textbf{ $h_{\tilde{\mu}}\left(S;\mathscr{\mathcal{A}}\right)=0$.}\\

\textbf{Proof:} Let $\tilde{\eta}$ be such an ergodic component,
then each $S$-factor corresponding to $\mathcal{A}_{k}$ is just
a periodic measure on a finite number of atoms and our $S$-system
with $\tilde{\eta}$ is isomorphic to the inverse limit of these $S$-factors.
In particular, our $S$-system with $\tilde{\eta}$ is isomorphic
to a Kronecker system (an ergodic rotation of a compact group) and
hence has $0$ entropy. $\blacksquare$\\

If $x,y\in X$ agree on the non-negative horizontal axis then $\nu_{n}\left(x\right)$
is just a translation of $\nu_{n}\left(y\right)$ by $\varphi\left(T^{-n}\left(x\right)\right)-\varphi\left(T^{-n}\left(y\right)\right)$.
Let us call a point $x\in X$ $symmetric$ if there exists $y\in X$
agreeing with it on the non-negative horizontal axis but disagreeing
with it on a coordinate $\left(-i_{0},0\right)$ for some $i_{0}>0$,
and $v_{n}\left(T^{m}\left(x\right)\right)=\nu_{n}\left(T^{m}\left(y\right)\right)$
for all $m\geq0$ and $n\in\mathbb{N}$ . This implies that $\nu_{n}\left(T^{m}\left(x\right)\right)$
is invariant under a translation by $\varphi\left(T^{-n+m}\left(x\right)\right)-\varphi\left(T^{-n+m}\left(y\right)\right)$.
\\

Next, we want to prove that $x$ is $\tilde{\mu}$-a.s. not symmetric
(Prop. 3.11).\\

\textbf{Lemma 3.9:} The set of symmetric points is $T$ and $S$ invariant
(and hence of measure either $0$ or $1$).\\

\textbf{Proof:} $T$ invariance is obvious. For $S$ invariance, let
$x,y$ be a pair of corresponding points as in the definition of a
symmetric point. We claim that $S\left(x\right)$ and $S\left(y\right)$
also form such a pair. Lemma 3.3 implies that if $i_{0}=\min\left\{ i\,:\,x\left(-i,0\right)\neq y\left(-i,0\right)\right\} $
then $S\left(x\right)\left(-i_{0},0\right)\neq S\left(y\right)\left(-i_{0},0\right)$.
So it suffices to prove that $v_{n}\left(T^{m}\left(S\left(x\right)\right)\right)=\nu_{n}\left(T^{m}\left(S\left(y\right)\right)\right)$,
and this will follow if we prove $a_{n}^{\left(qt\,\mod p^{n}\right)}\left(T^{m}S\left(x\right)\right)=a_{n}^{\left(qt\,\mod p^{n}\right)}\left(T^{m}S\left(y\right)\right)$
since $q$ is invertible in the ring $\mathbb{Z}/p^{n}\mathbb{Z}$
as

$\gcd\left(p,q\right)=1$ \footnote{This is the only place in the proof of Theorem 1.3 where the fact
$\gcd\left(p,q\right)=1$ is used not through Lemma 3.3.}.\\

Indeed, $a_{n}^{\left(qt\,\mod p^{n}\right)}\left(T^{m}S\left(x\right)\right)=a_{n}^{\left(t\right)}\left(T^{m}\left(x\right)\right)=a_{n}^{\left(t\right)}\left(T^{m}\left(y\right)\right)=a_{n}^{\left(qt\,\mod p^{n}\right)}\left(T^{m}S\left(y\right)\right)$
- this follows from Lemma 3.6 together with the fact that $T^{m}$
is a symmetry of the action of $S$. $\blacksquare$\\

\textbf{Lemma 3.10:} Given $r\in\mathbb{N}$, and integers $-r<c_{0},\dots,c_{n}<r$
with $c_{n}\neq0$ , if $\stackrel[i=0]{n}{\sum}\frac{c_{i}}{r^{i}}$
is equal to $\frac{u}{w}$ in least terms then $w\geq2^{n}$.\\

\textbf{Proof:} By induction on $n$. The case of $n=1$ is clear.
Continuing , if $r\frac{u}{w}=c_{0}+\stackrel[i=0]{n-1}{\sum}\frac{c_{i+1}}{r^{i}}$
equals $\frac{u'}{w'}$ in least terms then $w$ is a non-trivial
multiple of $w'$ since all prime divisors of $w$ divide $r$. Hence
$w\geq2w'\geq2^{n}$ - where the last inequality is by the induction
hypothesis. $\blacksquare$\\

\textbf{Proposition 3.11:} $\tilde{\mu}$-a.e. $x\in X$ is not symmetric.\\

\textbf{Proof:} By Lemma 3.9 and the ergodicity of $\tilde{\mu}$
to the $\mathbb{Z}^{2}$-action, if the propostion is false then $\tilde{\mu}$-a.e.
$x\in X$ is symmetric. We prove that this implies that $\mu$ is
Lebesgue measure $\lambda$ (in contradiction to our assumption in
the begining of this subsection).\\

We begin by showing that if $x\in X$ is a symmetric point then $\nu_{n}\left(x\right)\xrightarrow[n\rightarrow\infty]{}\lambda$
(weakly). So given such an $x$, there exists $x\neq y\in X$ agreeing
with it on the non-negative horizontal axis such that the measure
$\nu_{n}\left(x\right)$ is invariant under translation by

$\varphi\left(T^{-n}x\right)-\varphi\left(T^{-n}y\right)=\stackrel[i=1]{n}{\sum}\frac{y\left(-i,n-i\right)-x\left(-i,n-i\right)}{\left(pq\right)^{n-i+1}}$.
Let $i_{0}=\min\left\{ i\,:\,x\left(-i,0\right)\neq y\left(-i,0\right)\right\} $,
and assume $n\geq i_{0}$, then by Lemma 3.7, $\varphi\left(T^{-n}\left(x\right)\right)-\varphi\left(T^{-n}\left(y\right)\right)$
is a fraction with denomitor $\geq2^{n-i_{0}+1}$(in its least terms
representation). So the group of translations under which $\nu_{n}\left(x\right)$
is invariant is of order at least $2^{n-i_{0}+1}$ and thus contains
an element $\leq\frac{1}{2^{n-i_{0}+1}}$. This implies that

$\underset{\mathbb{T}}{\int}f\left(s\right)\,d\left(\nu_{n}\left(x\right)\right)\left(s\right)\xrightarrow[n\rightarrow\infty]{}\underset{\mathbb{T}}{\int}f\left(s\right)\,d\lambda\left(s\right)$
for every $f\in C\left(\mathbb{T}\right)$.\\

$\mu=\varphi_{*}\left(\underset{\mathbb{T}}{\int}\tilde{\mu}_{x}^{\stackrel[i=n]{\infty}{\vee}T^{-i}\left(\alpha\right)}\,d\tilde{\mu}\left(x\right)\right)=\underset{\mathbb{T}}{\int}\varphi_{*}\left(\tilde{\mu}_{x}^{\stackrel[i=n]{\infty}{\vee}T^{-i}\left(\alpha\right)}\right)\,d\tilde{\mu}\left(x\right)=\underset{\mathbb{T}}{\int}\left(\varphi\left(T^{-n}\left(x\right)\right)+\nu_{n}\left(x\right)\right)\,d\mu\left(x\right)$,
but $\varphi\left(T^{-n}\left(x\right)\right)+\nu_{n}\left(x\right)\xrightarrow[n\rightarrow\infty]{}\lambda$
(weakly) and so $\mu=\lambda$. $\blacksquare$\\

The following lemma is desired property (ii) of $\mathcal{A}$ and
thus concludes the proof of Theorem 1.3.\\

\textbf{Lemma 3.12:} $T\left(\alpha\right)\underset{\mod\tilde{\mu}}{\subseteq}\mathcal{A}\vee\left(\stackrel[i=0]{\infty}{\vee}T^{-i}\left(\alpha\right)\right)$.\\

\textbf{Proof:} The statement will follow if we prove that for $\tilde{\mu}$-a.e.
$x$, its atom with respect to

$\mathcal{A}\vee\left(\stackrel[i=0]{\infty}{\vee}T^{-i}\left(\alpha\right)\right)$
is contained in its atom with respect to $T\left(\alpha\right)$.
But if $x$ is not such and belongs to the full measure set on which
the functions $\nu_{n}\circ T^{m}$ are defined, it implies that there
exists $x\neq y\in X$ as in the definition of a symmetric point,
i.e. $x$ is symmetric, and by Lemma 3.11 we are done. $\blacksquare$\\

\section{Parry's Proof of Rudolph's Theorem}

This section is based on a proof of W. Parry to Rudolph's Theorem
\cite{key-19} (with some adjustments). The proof has similarities
to Rudolph's original proof (presented in section 3), but is still
different. Among the more superficial differences between the two
is that it is not formulated in the language of symbolic dynamics.
The only result from other sections we shall need here is Prop. 3.1.
\\

$\left(X,\mathcal{B},\mu\right)$ is said to be a \textit{standard
Borel probability space} if $X$ is a subset of a compact metric space
and $\mathcal{B}$ is its Borel $\sigma$-algebra.\\

\textbf{Lemma 4.1:} Assume $T$ is a surjective measure preserving
transformation on a Borel probability space $\left(X,\mathcal{B},\mu\right)$,
and that there exists a partition $\xi$ for which $\stackrel[i=0]{\infty}{\vee}T^{-i}\xi\stackrel{\mod\mu}{=}\mathcal{B}$.
Then every $T$-invariant sub-$\sigma$-algebra $\mathcal{A}$ - i.e.
$T^{-1}\mathcal{A}=\mathcal{A}$ - is contained in the Pinsker $\sigma$-algebra
of the system.\\

\textbf{Proof:} Note that

\[
h_{\mu}\left(T\,|\,\mathcal{A}\right)=H_{\mu}\left(\mathcal{\xi}\,|\,\stackrel[i=1]{\infty}{\vee}T^{-i}\mathcal{\xi}\vee\mathcal{A}\right)=H_{\mu}\left(\mathcal{\xi}\,|\,T^{-1}\mathcal{B}\vee T^{-1}\mathscr{\mathcal{A}}\right)=H_{\mu}\left(\mathcal{\xi}\,|\,T^{-1}\mathcal{B}\right)=h_{\mu}\left(T\right).
\]

The result follows by passing to the invertible extension (that is
the reason for the surjectivity requirement), and applying the Abramov-Rokhlin
formula. $\blacksquare$\\

\subsection{Invariance of Conditional Informations for Certain Commuting Maps}

Let $\left(X,\mathcal{B},\mu\right)$ be a Standard Borel probability
space, and $S:X\rightarrow X$ be a measurable map that preserves
$\mu$.\\

Any function measurable function $w:X\rightarrow\left[0,\infty\right]$
induces a measure $\mu_{w}$ defined by setting $d\mu_{w}=w\,d\mu$.
The measure $\mu_{w}$ can be pushed-foreward through $S$ to obtain
a new measure on $X$ which also is absolutely continuous with respect
to $\mu$, and let us denote its Radon-Nikodym derivative by $L_{S}w$
($L_{S}\left(\cdot\right)$ is called the transfer operator of $S$).
One may verify that
\[
\left(L_{S}w\right)\left(S\left(x\right)\right)=\mathbb{E}_{\mu}\left(w|S^{-1}\mathcal{B}\right)=\int w\left(y\right)\,d\mu_{x}^{S^{-1}\mathcal{B}}\left(y\right),
\]

where $\mu_{x}^{S^{-1}\mathcal{B}}$ is the conditional measure at
$x$ with respect to $S^{-1}\mathcal{B}$.

(Notice that in particular $L_{S}\left(w\circ S\right)\left(S\left(x\right)\right)=w\left(S\left(x\right)\right)$,
and because $S$ preserves $\mu$ this implies that $L_{S}\left(w\circ S\right)=w$
- but we shall not make use of this fact.)\\

We shall need the following relation between the transfer operator
and the Koopman operator of $S$.\\

\textbf{Proposition 4.2:} Let $\left(X,\mathcal{B},\mu\right)$ be
a Standard Borel probability space, $S:X\rightarrow X$ a measurable
map that preserves $\mu$ and $w\in L_{\mu}^{1}\left(X\right)$. If
$L_{S}\left(w\right)=w$ then $w\circ S=w$. \footnote{This proposition is well known (see for example exercise 10.19 in
\cite{key-9}). However, the author learned the specific proof presented
here from ChatGPT, and neither he or it could find a reference for
this exact proof. In addition, ChatGPT had also pointed out how to
deduce from this proposition its truth for a measurable non-negative
$w$ that may not be integrable, however we do not need here this
strengthening.}\\

\textbf{Proof:} Notice that
\[
w\circ S^{n}=L_{S}^{n}\left(w\right)\circ S^{n}=L_{S^{n}}\left(w\right)\circ S^{n}=\mathbb{E}_{\mu}\left(w|S^{-n}\mathcal{B}\right).
\]

So by the reverse martingale convergence theorem, the sequence $w\circ S^{n}$
converges in $L_{\mu}^{1}\left(X\right)$. In particular,
\[
\left\Vert w\circ S^{n+1}-w\circ S^{n}\right\Vert _{1}\xrightarrow[n\rightarrow\infty]{}0.
\]

However, since composition by $S$ is an isometry ($S$ is measure
preserving), we have that
\[
\left\Vert w\circ S-w\right\Vert _{1}\xrightarrow[n\rightarrow\infty]{}0,
\]

i.e. that
\[
w\circ S=w.\,\,\,\,\,\,\,\,\,\,\,\,\,\,\,\,\,\,\blacksquare
\]

If $S$ is countable-to-one and we define $f\left(x\right)=I_{\mu}\left(\mathcal{B}\,|\,S^{-1}\mathcal{B}\right)\left(x\right)$
(the conditional information function), then $L_{S}w\left(x\right)=\underset{y\in S^{-1}x}{\sum}e^{-f\left(y\right)}w\left(y\right)$
(under the convention $0\cdot\infty=0$).\\

Assuming this and throwing into the game another measurable map $T:X\rightarrow X$
satisfying for some $M>0$ that $\left|T^{-1}x\right|\leq M$ for
all $x\in X$, preserves $\mu$ and commutes with $S$ we define $g\left(x\right)=I\left(\mathcal{B}\,|\,T^{-1}\mathcal{B}\right)\left(x\right)$,
and so $L_{T}w\left(x\right)=\underset{y\in T^{-1}x}{\sum}e^{-g\left(y\right)}w\left(y\right)$.
The commutation relation implies the important identity

\[
\left(**\right)\,\,\,\,f\left(Tx\right)-f\left(x\right)=g\left(Sx\right)-g\left(x\right),
\]
 proved by exponentiating both sides of the equation $-f\left(x\right)-g\left(Sx\right)=-g\left(x\right)-f\left(Tx\right)$,
which then reads 
\[
\mu_{x}^{S^{-1}\mathcal{B}}\left(\left\{ x\right\} \right)\cdot\mu_{Sx}^{T^{-1}\mathcal{B}}\left(\left\{ Sx\right\} \right)=\mu_{x}^{T^{-1}\mathcal{B}}\left(\left\{ x\right\} \right)\cdot\mu_{Tx}^{S^{-1}\mathcal{B}}\left(\left\{ Tx\right\} \right),
\]
i.e. $\mu_{x}^{\left(TS\right)^{-1}\mathcal{B}}\left(\left\{ x\right\} \right)=\mu_{x}^{\left(ST\right)^{-1}\mathcal{B}}\left(\left\{ x\right\} \right)$.\\

We now assume further that for $\mu$-almost-every $x$ the map $T:S^{-1}x\rightarrow S^{-1}Tx$
is bijective, and show that this implies that\uline{ \mbox{$f\circ T=f$}
and \mbox{$g\circ S=g$}}.\\
\\
 $\left(L_{S}e^{g\left(\cdot\right)}\right)\left(x\right)=\underset{y\in S^{-1}\left(x\right)}{\sum}e^{-f\left(y\right)+g\left(y\right)}$,
so by identity $\left(**\right)$ this is equal to 
\[
\underset{y\in S^{-1}x}{\sum}e^{g\left(S\left(y\right)\right)-f\left(T\left(y\right)\right)}=e^{g\left(x\right)}\cdot\underset{y\in S^{-1}x}{\sum}e^{-f\left(T\left(y\right)\right)}=e^{g\left(x\right)}
\]

(where in the last equality we applied the bijectivity assumption).\\

Hence $e^{g\left(x\right)}$ is fixed by $L_{S}$. Since
\[
\int e^{g\left(x\right)}d\mu\left(x\right)=\int\frac{1}{\mu_{x}^{T^{-1}\mathcal{B}}\left(\left\{ x\right\} \right)}d\mu\left(x\right)=\int\int\frac{1}{\mu_{y}^{T^{-1}\mathcal{B}}\left(\left\{ y\right\} \right)}d\mu_{x}^{T^{-1}\mathcal{B}}\left(y\right)d\mu\left(x\right)\leq\int Md\mu\left(x\right)=M,
\]

the function $e^{g\left(x\right)}$ belongs to $L_{\mu}^{1}\left(X\right)$.
Hence by Prop. 4.2 $e^{g\left(S\left(x\right)\right)}=e^{g\left(x\right)}$.
So $g\left(S\left(x\right)\right)=g\left(x\right)$, and by $\left(**\right)$
also $f\left(T\left(x\right)\right)=f\left(x\right)$.\\

\subsection{The Proof}

Sticking to Parry's original notation, let $p,q>1$ be relatively
prime integers and $S,T:\mathbb{T}\rightarrow\mathbb{T}$ be the multiplications
by $p,q$ respectively.\\

Let $\mu\in\text{\ensuremath{M\left(\mathbb{T}\right)}}$ be invariant
under both $S$ and $T$ and ergodic under their joint action. We
assume further it is of positive entropy with respect to the action
of each and need to prove that it is Lebesgue measure.\\

The fact that $p$ and $q$ are relatively prime assures that $T:S^{-1}x\rightarrow S^{-1}Tx$
is a bijection for every $x\in\mathbb{T}$. Hence all conclusions
of section 4.1 apply in this case. In particular, $f\circ T=f$.\\

For $\mu$-almost-every $x\in\mathbb{T}$ and $n\in\mathbb{N}$ we
define $d\left(x,n\right)\in\ensuremath{M\left(\mathbb{T}\right)}$
supported on $0,\frac{1}{p^{n}},\frac{2}{p^{n}},\dots,\frac{p^{n}-1}{p^{n}}$
where the mass of each $\frac{i}{p^{n}}$ is $\mu_{x}^{S^{-n}\mathcal{B}}\left(\left\{ x+\frac{i}{p^{n}}\right\} \right)$.

Notice that hence $d\left(x,n\right)$ determines $d\left(x,n-1\right)$.
Also, if we define

$f^{n}=f+f\circ S+\dots+f\circ S^{n-1}$ then $\mu_{x}^{S^{-n}\mathcal{B}}\left(\left\{ x+\frac{i}{p^{n}}\right\} \right)=e^{-f^{n}\left(x+\frac{i}{p^{n}}\right)}$.
\\

We denote by $\mathscr{\mathcal{H}}_{n}\subseteq\mathscr{\mathcal{B}}$
the smallest $\sigma$-algebra on $X$ for which the function $d\left(\cdot,n\right):\mathbb{T}\rightarrow M\left(\mathbb{T}\right)$
is measurable (so $\mathscr{\mathcal{H}}_{1}\subseteq\mathscr{\mathcal{H}}_{2}\subseteq\mathscr{\mathcal{H}}_{3}\subseteq\dots$).
Since $f\circ T=f$ and $S,T$ commute then $f^{n}\circ T=f^{n}$.
\\

Applying this we note that
\[
f^{n}\left(x+\frac{i}{p^{n}}\right)=f^{n}\left(T\left(x+\frac{i}{p^{n}}\right)\right)=f^{n}\left(Tx+\frac{qi}{p^{n}}\right),
\]

namely that $\mu_{x}^{S^{-n}\mathcal{B}}\left(\left\{ x+\frac{i}{p^{n}}\right\} \right)=\mu_{Tx}^{S^{-n}\mathcal{B}}\left(\left\{ Tx+\frac{qi}{p^{n}}\right\} \right)$.
Thus $T^{-1}\mathscr{\mathcal{H}}_{n}=\mathscr{\mathcal{H}}_{n}$.
Defining

$\mathscr{\mathcal{H}}=\stackrel[n=1]{\infty}{\vee}\mathscr{\mathcal{H}}_{n}$
we obtain that $\mathscr{\mathcal{H}}$ is $T$-invariant (i.e. $T^{-1}\mathscr{\mathcal{H}}=\mathscr{\mathcal{H}}$).
By Lemma 4.1, this means that $\mathscr{\mathcal{H}}$ is contained
in the Pinsker $\sigma$-algebra of $T$. \\

\textbf{Lemma 4.3:} The Pinsker $\sigma$-algebras of $S$ and $T$
are equal.\\

\textbf{Proof:} Denote by $\tilde{\mu}$ the measure in the invertible
extension of our system relative to the joint action of both $S$
and $T$, each of the corresponding maps there $\tilde{S},\tilde{T}$
is an automorphism of the system defined by the other and hence the
Pinsker $\sigma$-algebra of each of the maps there is invariant relative
to the other. Denote these Pinsker $\sigma$-algebras by $\mathcal{P}_{\tilde{S}},\mathcal{P}_{\tilde{T}}$.
We claim $\mathcal{P}_{\tilde{T}}=\mathcal{P}_{\tilde{S}}$. By a
similar proof to Prop. 3.1 $h_{\tilde{\mu}}\left(\tilde{T}\,|\,\mathcal{P}_{\tilde{T}}\right)=\frac{\log p}{\log q}h_{\tilde{\mu}}\left(\tilde{S}\,|\,\mathcal{P}_{\tilde{T}}\right)$
but $h_{\tilde{\mu}}\left(\tilde{T}\,|\,\mathcal{P}_{\tilde{T}}\right)=h_{\tilde{\mu}}\left(\tilde{T}\right)=\frac{\log p}{\log q}h_{\tilde{\mu}}\left(\tilde{S}\right)$
(the latter equality is implied by Prop. 3.1) and hence $h_{\tilde{\mu}}\left(\tilde{S}\,|\,\mathcal{P}_{\tilde{T}}\right)=h_{\tilde{\mu}}\left(\tilde{S}\right)$.
By the Abaramov-Rokhlin formula $\mathcal{P}_{\tilde{T}}\subseteq\mathcal{P}_{\tilde{S}}$.
In the same manner one shows the reverse containment.\\

Now, if a set $A$ belongs to the Pinsker $\sigma$-algebra down at
the circle of $T$ (resp. $S$) then its inverse image up at the joint
action invertible extension belongs to the Pinsker $\sigma$-algebra
of $T$ (resp. $S$) there and hence of $S$ (resp. $T$) there, and
this means that $A$ belongs to the Pinsker $\sigma$-algebra of $S$
(resp. $T$) down at the circle. \textbf{$\blacksquare$}\\

Hence, by Lemma 4.3, $\mathcal{H}\subseteq\stackrel[n=1]{\infty}{\cap}S^{-n}\mathcal{B}$
(modulo $\mu$). This last containment is the key to the proof.\\

There exists a set $N\subseteq\mathbb{T}$ of measure $0$ such that
the conditional measures $\mu_{x}^{S^{-n}\mathscr{\mathcal{B}}}$
are defined on $\mathbb{T}\setminus N$ for all $n\in\mathbb{N}$,
and $\mu_{x}^{S^{-n}\mathscr{\mathcal{B}}}\left(\left[x\right]_{S^{-n}\mathscr{\mathcal{B}}}\setminus N\right)=1$,
$\left[x\right]_{S^{-n}\mathscr{\mathcal{B}}}\setminus N\subseteq\left[x\right]_{\mathscr{\mathcal{H}}}$
for all $x\in\mathbb{T}\setminus N$ (where the square brackets denote
an atom of the $\sigma$-algebra appearing in the subscript) - here
we use the fact that $\mathcal{H}\subseteq\stackrel[n=1]{\infty}{\cap}S^{-n}\mathcal{B}$.\\

Thus $\mu_{x}^{S^{-n}\mathscr{\mathcal{B}}}$ is the uniform measure
on $\left[x\right]_{S^{-n}\mathscr{\mathcal{B}}}\setminus N$ (each
element with probability $e^{-f^{n}\left(x\right)}$) for every $x\in\mathbb{T}\setminus N$.
Moreover, given any $x\in\mathbb{\mathbb{T}\setminus N}$ and any
$y\in\left[x\right]_{S^{-1}\mathscr{\mathcal{B}}}\setminus N$ the
equality $d\left(y,n\right)=d\left(y+g,n\right)$ holds for any $g$
which is a difference of two elements of the set $\left[x\right]_{S^{-n}\mathscr{\mathcal{B}}}\setminus N$.
Therefore for every $x\in\mathbb{T}\setminus N$ there exists a subgroup
$G_{x}^{n}$ of $\left\{ 0,\frac{1}{p^{n}},\frac{2}{p^{n}},\dots,\frac{p^{n}-1}{p^{n}}\right\} $
such that $\left[x\right]_{S^{-n}\mathscr{\mathcal{B}}}\setminus N$
is a coset of $G_{x}^{n}$.\\

$h_{\mu}\left(S\right)>0$ implies that the set $E\subseteq\mathbb{T}\setminus N$
composed of the points $x$ such that $\left|\left[x\right]_{S^{-1}\mathscr{\mathcal{B}}}\setminus N\right|>1$
is of positive measure, so for them $G_{x}^{n}\neq\left\{ 0\right\} $
for every $n$. But $X\setminus E\subseteq T^{-1}\left(X\setminus E\right)$
(modulo $\mu$) and, since the joint action of $S$ and $T$ is ergodic,
the Pointwise Ergodic Theorem implies that $S^{n}x\in E$ infinitely
often for $x\in\mathbb{T}$ almost-surely, and thus $f^{n}\left(x\right)\underset{n\rightarrow\infty}{\longrightarrow}\infty$
almost-surely. But $e^{f^{n}\left(x\right)}=\left|G_{x}^{n}\right|$
and thus $\left|G_{x}^{n}\right|\underset{n\rightarrow\infty}{\longrightarrow}\infty$
almost-surely.\\

So $\mathbb{E}_{\mu}\left(e^{2\pi kxi}\,|\,S^{-n}\mathscr{\mathcal{B}}\right)\left(x\right)=\frac{\underset{g\in G_{x}^{n}}{\sum}e^{2\pi k\left(x+g\right)i}}{\left|G_{x}^{n}\right|}$
equals $0$ for $n$ suffiecients large and this guarantees $\mathbb{E}_{\mu}\left(e^{2\pi kxi}\right)=0$
for every integer $k\neq0$. This means $\mu$ is Lebesgue measure.\\

\section{Host's Theorem}

The main source for the discussion of Host's Theorem presented here
is \cite{key-18} by D. Meiri (it is both a mathematical and a linguistical
strenghtening of B. Host's original paper which is in French \cite{key-15}
- for predecessors see \cite{key-10,key-21}).\\

\textbf{Definition:} For an integer $p>1$, a sequence $c_{k}\in\mathbb{N}$
is called a $p$-Host sequence if for every $\mu\in\text{\ensuremath{M\left(\mathbb{T}\right)}}$
invariant under multiplication by $p$, ergodic and with entropy $h>0$,
the sequence $c_{k}x$ equidistributes $\mu$-a.s..\\

Note that the $\mu$-a.s. equidistribution of $c_{k}x$ implies that
every $f\in C\left(\mathbb{T}\right)$ satisfies

$\frac{1}{N}\sum_{k=0}^{N-1}\int f\left(x\right)\,d\left(c_{k}\mu\right)=\frac{1}{N}\sum_{k=0}^{N-1}\int f\left(c_{k}x\right)\,d\mu=\int\frac{1}{N}\sum_{k=0}^{N-1}f\left(c_{k}x\right)\,d\mu\xrightarrow[N\rightarrow\infty]{}\int f\left(x\right)\,d\lambda$
(where $\lambda$ is Lebesgue measure), i.e. the sequence $c_{k}\mu$
equidistributes.\\

\textbf{Theorem 5.1 (Host's Theorem):} Let $p,q>1$ be relatively
prime integers, then $q^{k}$ is a $p$-Host sequence. \footnote{The stronger version - for $p,q$ multiplicatively independent, was
proved a few years ago by Hochman and Shmerkin \cite{key-13}. }\\

Before presenting the proof, let us see why Host's Theorem implies
Rudolph's Theorem (and thus serves as an independent proof of the
latter). Given an atomless $\mu\in\text{\ensuremath{M\left(\mathbb{T}\right)}}$
invariant and ergodic with respect to the action of the multiplicative
semi-group $\left\langle p,q\right\rangle $ and $h_{\mu}\left(p\right)>0$,
consider the ergodic decomposition of $\mu$ with respect to the times
$p$ transformation: $\mu=\int_{\mathbb{T}}\mu_{x}^{\varepsilon}\,d\mu\left(x\right)$
($\mu_{x}^{\varepsilon}$ is the conditional measure at $x$ with
respect to $\varepsilon$ - the $\sigma$-algebra of invariant sets).
By averaging with respect to iterates of the times $q$ transformation
we obtain $\mu=\frac{1}{N}\Sigma_{k=0}^{N-1}q^{k}\mu=\int_{\mathbb{T}}\frac{1}{N}\Sigma_{k=0}^{N-1}q^{k}\mu_{x}^{\varepsilon}\,d\mu\left(x\right)$.
Now $h_{\mu}\left(p\right)=\int_{\mathbb{T}}h_{\mu_{x}^{\varepsilon}}\left(p\right)\,d\mu\left(x\right)$
and thus there exists $A\subseteq\mathbb{T}$ with $\mu\left(A\right)>0$
for which $h_{\mu_{x}^{\varepsilon}}\left(p\right)>0$ for every $x\in A$.
By Host's Theorem $\int_{A}\frac{1}{N}\Sigma_{k=0}^{N-1}q^{k}\mu_{x}^{\varepsilon}\,d\mu\left(x\right)$
converges to $\mu\left(A\right)\lambda$, and so $\mu\gg\lambda$
(if $\mu\left(B\right)=0$ then $\mu_{x}^{\varepsilon}\left(B\right)=0$
$\mu$-a.s.), but $\mu$ is ergodic with respect to the semi-group
action and hence $\mu=\lambda$.\\

We need one number theoretical fact which is an immediate conclusion
from Lemma 2.3 (this is the only place in this section where the comprimality
of $p,q$ is exploited).\\

\textbf{Proposition 5.2: }For every $0\neq a\in\mathbb{Z}$ there
exists some $M>0$ for which 

$\#\left\{ k\,:\,0\leq k<p^{n},\,aq^{k}\equiv t\,\left(\mod p^{n}\right)\right\} <M$
for all integers $n>0$ and $t$.\\

We now turn to prove Host's Theorem. Let $\mu$ be a measure such
as in the definition of a $p$-Host sequence. In order to prove the
theorem it is enough to show that $g_{N}\left(x\right)=\frac{1}{N}\sum_{k=0}^{N}e\left(aq^{k}x\right)\xrightarrow[N\rightarrow\infty]{}0$
$\mu$-a.s. for every $0\neq a\in\mathbb{Z}$, where $e\left(x\right)=e^{2\pi xi}$
for $x\in\mathbb{T}=\mathbb{R}/\mathbb{Z}$. $\int\left|g_{N}\left(x\right)\right|^{2}\,d\mu$
is bounded by $1$, but this does not solve the problem. As we shall
now see, Host observed that if one considers a sum of translates $\omega_{n\left(N\right)}=\sum_{j=0}^{p^{n\left(N\right)}-1}\delta_{\frac{j}{p^{n\left(N\right)}}}*\mu$
instead of $\mu$, the integral of $\left|g_{N}\left(x\right)\right|^{2}$
by this measure (which is of total mass $p^{n}$) is surprisingly
still bounded uniformly in $N$ by a constant due to cancelations.\\

We choose $n\left(N\right)$ to be the natural number satisfying $p^{n-1}\leq N<p^{n}$.
Evaluating the expression\\

$\int\left|g_{N}\left(x\right)\right|^{2}\,d\omega_{n}=\int\sum_{j=0}^{p^{n}-1}\left|g_{N}\left(x+\frac{j}{p^{n}}\right)\right|^{2}\,d\mu=\frac{1}{N^{2}}\sum_{k,l=0}^{N-1}\sum_{j=0}^{p^{n}-1}e\left(\left(aq^{k}-aq^{l}\right)\frac{j}{p^{n}}\right)\cdot\left(\int e\left(a\left(q^{k}-q^{l}\right)x\right)\,d\mu\right)$\\

$\leq\frac{1}{N^{2}}\sum_{k,l=0}^{N-1}\left|\sum_{j=0}^{p^{n}-1}e\left(\left(aq^{k}-aq^{l}\right)\frac{j}{p^{n}}\right)\right|\cdot\left|\int e\left(a\left(q^{k}-q^{l}\right)x\right)\,d\mu\right|\leq$\\

$\leq\frac{1}{N^{2}}\sum_{k,l=0}^{N-1}\left|\sum_{j=0}^{p^{n}-1}e\left(\left(aq^{k}-aq^{l}\right)\frac{j}{p^{n}}\right)\right|$.\\

The summation over $j$ vanishes if $aq^{k}\not\equiv aq^{l}\,\left(\mod\,p^{n}\right)$
and is $p^{n}$ otherwise, hence\\

$\int\left|g_{N}\left(x\right)\right|^{2}\,d\omega_{n}\leq\frac{p^{n}}{N^{2}}\cdot\#\left\{ \left(k,l\right)\in\left\{ 0,...,N-1\right\} ^{2}:\,aq^{k}\equiv aq^{l}\,\left(\mod p^{n}\right)\right\} $\\

$=\frac{p^{n}}{N^{2}}\sum_{t=0}^{p^{n}-1}\left(\#\left\{ 0\leq k\leq N-1\,:\,aq^{k}\equiv t\,\left(\mod p^{n}\right)\right\} \right)^{2}$.\\

Proposition 5.2 implies that $\int\left|g_{N}\left(x\right)\right|^{2}\,d\omega_{n}\leq\frac{p^{2n}M^{2}}{N^{2}}\leq\frac{N^{2}p^{2}M^{2}}{N^{2}}=p^{2}M^{2}$
for all $N$.\\

On the other side $\int_{\frac{d\mu}{d\omega_{n}}>0}\frac{\left|g_{N}\left(x\right)\right|^{2}}{\frac{d\mu}{d\omega_{n}}\left(x\right)}\,d\mu=\int_{\frac{d\mu}{d\omega_{n}}>0}\left|g_{N}\left(x\right)\right|^{2}\frac{d\omega_{n}}{d\mu}\left(x\right)\,d\mu\leq\int\left|g_{N}\left(x\right)\right|^{2}\,d\omega_{n}$,
so we conclude $\int\frac{\left|g_{N}\left(x\right)\right|^{2}}{\frac{d\mu}{d\omega_{n}}\left(x\right)}\,d\mu\leq p^{2}M$
($\frac{d\mu}{d\omega_{n}}\left(x\right)>0$ $\mu$-a.s. hence, in
particular, the integrand is well-defined). This already looks interesting.\\

Let us now investigate the Radon-Nikodym derivative $\frac{d\mu}{d\omega_{n}}$.
Denote by $\alpha$ the partition $\left\{ \left[\frac{j}{p},\frac{j+1}{p}\right)\right\} _{j=0}^{p-1}$,
and by $\alpha_{k}^{k+l}$ denote $\sigma_{p}^{-k}\left(\alpha\right)\vee...\vee\sigma_{p}^{-\left(k+l\right)}\left(\alpha\right)$,
where $\sigma_{p}$ is the transformation of multiplication by $p$.\\

\textbf{Proposition 5.3:} $-\log\frac{d\mu}{d\omega_{n}}=I_{\mu}\left(\alpha_{0}^{n-1}|\alpha_{n}^{\infty}\right)$
(the conditional information function) $\mu$-a.s..\\

\textbf{Proof:} It is sufficient to prove that $\mathbb{\mathbb{E}}_{\mu}\left(f|\alpha_{n}^{\infty}\right)\left(x\right)=\sum_{j=0}^{p^{n}-1}f\left(x+\frac{j}{p^{n}}\right)\frac{d\mu}{d\omega_{n}}\left(x+\frac{j}{p^{n}}\right)$
$\mu$-a.s. for any $f\in L^{1}\left(\mathbb{T},\mu\right)$ (then
for every element $S$ of the partition $\alpha_{0}^{n}$ one takes
$f=1_{S}$). The function on the right hand side is indeed $\alpha_{n}^{\infty}$-measurable,
and for every set $A\in\alpha_{n}^{\infty}$\\

$\int_{A}f\left(x\right)\,d\mu\left(x\right)=\int_{A}f\left(x\right)\frac{d\mu}{d\omega_{n}}\left(x\right)\,d\omega_{n}=\sum_{j=0}^{p^{n}-1}\int_{A}f\left(x+\frac{j}{p^{n}}\right)\frac{d\mu}{d\omega_{n}}\left(x+\frac{j}{p^{n}}\right)\,d\mu\left(x\right)$.
$\blacksquare$\\

The following proposition is reminiscent of the Shannon-Mcmillan-Breiman
Theorem (although it is much easier to prove).\\

\textbf{Proposition 5.4:} $-\frac{1}{n}\log\frac{d\mu}{d\omega_{n}}\xrightarrow[n\rightarrow\infty]{}h$
$\mu$-a.s..\\

\textbf{Proof:} The proof is a simple implimentation of the Pointwise
Ergodic Theorem:\\

$-\frac{1}{n}\log\frac{d\mu}{d\omega_{n}}=\frac{1}{n}I_{\mu}\left(\alpha_{0}^{n-1}|\alpha_{n}^{\infty}\right)=\frac{1}{n}\left(I_{\mu}\left(\alpha|\alpha_{1}^{\infty}\right)+I_{\mu}\left(\alpha_{1}|\alpha_{2}^{\infty}\right)+...+I_{\mu}\left(\alpha_{n-1}|\alpha_{n}^{\infty}\right)\right)$\\

$=\frac{1}{n}\left(I_{\mu}\left(\alpha|\alpha_{1}^{\infty}\right)+I_{\mu}\left(\alpha|\alpha_{1}^{\infty}\right)\circ\sigma_{p}+...+I_{\mu}\left(\alpha|\alpha_{1}^{\infty}\right)\circ\sigma_{p}^{n-1}\right)$.
$\blacksquare$\\

Summing up what we know up to now: $\left(1\right)$ $\left(\frac{d\mu}{d\omega_{n}}\right)^{\frac{1}{n}}\xrightarrow[n\rightarrow\infty]{}e^{-h}$
$\mu$-a.s. and $\left(2\right)$ $\int\frac{\left|g_{N}\left(x\right)\right|^{2}}{\frac{d\mu}{d\omega_{n}}\left(x\right)}\,d\mu$
is bounded as a sequence in $N$ ($n$ is determined by the condition
$p^{n-1}\leq N<p^{n}$). Given any $C>1$, $\left(2\right)$ implies
that\\

$\int\sum_{N=0}^{\infty}\frac{\left|g_{\left\lfloor C^{N}\right\rfloor }\left(x\right)\right|^{2}}{e^{\frac{nh}{3}}\frac{d\mu}{d\omega_{n}}\left(x\right)}\,d\mu=\sum_{N=0}^{\infty}\int\frac{\left|g_{\left\lfloor C^{N}\right\rfloor }\left(x\right)\right|^{2}}{e^{\frac{nh}{3}}\frac{d\mu}{d\omega_{n}}\left(x\right)}\,d\mu<\infty$
and so $\frac{\left|g_{\left\lfloor C^{N}\right\rfloor }\left(x\right)\right|^{2}}{e^{\frac{nh}{3}}\frac{d\mu}{d\omega_{n}}\left(x\right)}\xrightarrow[N\rightarrow\infty]{}0$
$\mu$-a.s. - where here we we define $n$ by $p^{n-1}\leq\left\lfloor C^{N}\right\rfloor <p^{n}$.
But $\left(1\right)$ implies that for $\mu$-almost every $x$ the
inequality $\frac{d\mu}{d\omega_{n}}\left(x\right)\leq e^{-\frac{nh}{2}}\left(x\right)$
is satisfied for $n$ large enough, and thus for $\mu$-almost every
$x$ the expression $e^{\frac{nh}{3}}\frac{d\mu}{d\omega_{n}}\left(x\right)$
is bounded. This means $\left|g_{\left\lfloor C^{N}\right\rfloor }\left(x\right)\right|\xrightarrow[N\rightarrow\infty]{}0$
for $\mu$-almost every $x$. It is only in these last steps that
we use the positivity of $h$.\\

I learned from E. Lindenstrauss the following simple and useful lemma
from which one can easily deduce that $\left|g_{N}\left(x\right)\right|\xrightarrow[N\rightarrow\infty]{}0$
$\mu$-a.s. and this finishes the proof of Host's Theorem.\\

\textbf{Lemma 5.5:} Given any real number $C>1$, if $\frac{1}{\left\lfloor C^{n}\right\rfloor }\sum_{k=0}^{\left\lfloor C^{n}\right\rfloor -1}e\left(x_{k}\right)\xrightarrow[n\rightarrow\infty]{}0$
for a sequence $x_{k}\in\mathbb{T}$ then $\limsup_{N\rightarrow\infty}\left|\frac{1}{N}\sum_{k=0}^{N-1}e\left(x_{k}\right)\right|\leq1-\frac{1}{C}$.\\

\textbf{Proof:} Let $N$ be a natural number, and take $n$ for which
$\left\lfloor C^{n-1}\right\rfloor \leq N<\left\lfloor C^{n}\right\rfloor $.\\

Doing the math $\left|\frac{1}{\left\lfloor C^{n}\right\rfloor }\sum_{k=0}^{\left\lfloor C^{n}\right\rfloor -1}e\left(x_{k}\right)-\frac{1}{N}\sum_{k=0}^{N-1}e\left(x_{k}\right)\right|$\\

$=\left|\frac{1}{\left\lfloor C^{n}\right\rfloor }\sum_{k=0}^{\left\lfloor C^{n-1}\right\rfloor -1}e\left(x_{k}\right)-\frac{1}{N}\sum_{k=0}^{\left\lfloor C^{n-1}\right\rfloor -1}e\left(x_{k}\right)+\left(\frac{1}{\left\lfloor C^{n}\right\rfloor }-\frac{1}{N}\right)\sum_{k=\left\lfloor C^{n-1}\right\rfloor }^{N-1}e\left(x_{k}\right)+\frac{1}{\left\lfloor C^{n}\right\rfloor }\sum_{k=N}^{\left\lfloor C^{n}\right\rfloor -1}e\left(x_{k}\right)\right|$.\\

Each of the first two summands converges to $0$ when $N\rightarrow\infty$.
The sum of the two others is less in absolute value than\\

$\left(\frac{1}{\left\lfloor C^{n}\right\rfloor }-\frac{1}{N}\right)\left(\left\lfloor C^{n}\right\rfloor -\left\lfloor C^{n-1}\right\rfloor \right)+\frac{\left\lfloor C^{n}\right\rfloor -\left\lfloor C^{n-1}\right\rfloor }{\left\lfloor C^{n}\right\rfloor }\leq\left(\frac{\left\lfloor C^{n-1}\right\rfloor }{\left\lfloor C^{n}\right\rfloor }-\frac{1}{C}\right)\left(\frac{\left\lfloor C^{n}\right\rfloor }{\left\lfloor C^{n-1}\right\rfloor }-1\right)+1-\frac{\left\lfloor C^{n-1}\right\rfloor }{\left\lfloor C^{n}\right\rfloor }$\\

$=2-2\frac{\left\lfloor C^{n-1}\right\rfloor }{\left\lfloor C^{n}\right\rfloor }-\frac{\left\lfloor C^{n}\right\rfloor }{C\left\lfloor C^{n-1}\right\rfloor }+\frac{1}{C}$
which converges to $1-\text{\ensuremath{\frac{1}{C}}}.$ $\blacksquare$\\

Einstein Institute of Mathematics, Edmond J. Safra campus, The Hebrew
University of Jerusalem, Israel.\\

matan.tal@mail.huji.ac.il
\end{document}